\documentclass{amsart}
 
 \usepackage{amsmath, amscd, graphicx, amsbsy, pstricks}
\usepackage{epsfig}
\usepackage[cmtip,all]{xy}

\newcommand{\MathVector}[1]{\boldsymbol{#1}}
\newcommand{\intersect}[1]{\bigl< #1 \bigr>}
 
\newcommand{\hFunction}{\mathfrak{h}}
\newcommand{\FiniteCount}[1]{\#\left| #1 \right|}
\newcommand{\DeformedVirasoro}{\hat{V}}
\newcommand{\Virasoro}{V}

\newcommand{\Moduli}{\mathcal{M}}
\newcommand{\CompactModuli}{\overline{\mathcal{M}}}
\newcommand{\Teichmuller}{\mathcal{T}}

\newcommand{\Integers}{\mathbb{Z}}
\newcommand{\Reals}{\mathbb{R}}
\DeclareMathOperator{\MCG}{Mod}
\DeclareMathOperator{\Stab}{Stab}
\DeclareMathOperator{\Vol}{Vol}
\DeclareMathOperator{\Curvature}{Curv}
\DeclareMathOperator{\trace}{tr}

\newtheorem{theorem}{Theorem}[section]
\newtheorem{proposition}[theorem]{Proposition}

\numberwithin{equation}{section}

\allowdisplaybreaks

\begin{document}
\title[Mirzakhani's Recursion and Virasoro Constraints]{Mirzakhani's Recursion Relations, Virasoro Constraints and the KdV Hierarchy}

\author[Motohico Mulase]{Motohico Mulase$^1$}
\address{Department of Mathematics, University of California, Davis, CA 95616-8633}
\email{mulase@math.ucdavis.edu}
\thanks{$^1$Research supported by NSF grant DMS-0406077 and UC Davis.}

\author[Brad Safnuk]{Brad Safnuk$^2$}
\address{Department of Mathematics, University of California, Davis, CA 95616-8633}
\email{safnuk@math.ucdavis.edu}
\thanks{$^2$Research supported by NSF grant DMS-0406077 and a UC Davis Dissertation Year Fellowship.}



\begin{abstract}
We present in this paper a differential version of 
Mirzakhani's recursion relation for the Weil-Peters\-son volumes of the moduli spaces of bordered Riemann surfaces. We discover that
the differential relation, which is equivalent to the original integral
formula of Mirzakhani,
 is a Virasoro constraint
condition on a generating function for these volumes. 
 We also show that the generating function for $\psi$ and $\kappa_1$ intersections on $\CompactModuli_{g,n}$ is a 1-parameter solution to the KdV hierarchy.  It recovers the Witten-Kontsevich generating function
 when the parameter is set to be $0$. 
\end{abstract}
\maketitle


 \section{Introduction}
 
 In her striking series of papers \cite{art:MirzakhaniWP, art:MirzakhaniIT},
 Mirzakhani obtained a beautiful recursion formula  for the 
 Weil-Petersson volume of the moduli spaces of 
 bordered Riemann surfaces. Her recursion relation is an 
 integral formula involving a kernel function that
 appears in the work of McShane  \cite{art:McShane}
 on hyperbolic geometry of 
 surfaces.
 We have discovered that the differential version of 
 the Mirzakhani recursion formula, which is
 \emph{equivalent} to the original integral form, is indeed
 a Virasoro constraint condition imposed on a generating 
 function of these volumes. 
 
 Mirzakhani proves in \cite{art:MirzakhaniIT}
 that her recursion relation reduces to the Virasoro constraint 
 condition as the length parameters of the boundary components
 of Riemann surfaces go to infinity, and moreover, it recovers the 
 celebrated Witten-Kontsevich theorem of intersection numbers
 of tautological classes on the moduli spaces of stable algebraic curves.
 Our result reveals that the Virasoro structure exists essentially in 
 the Mirzakhani theory, and that it is not the consequence of the large 
 boundary limit. 
 
 The Virasoro constraint formulas for the generating functions
 of Gromov-Witten invariants of various target manifolds have
 been extensively studied in recent years \cite{art:DubrovinZhang, art:Getzler, art:Givental, art:OK-GW}.
 Although Mirzakhani's hyperbolic method does not immediately
 apply to these cases with higher dimensional target spaces, the
 Virasoro structure we identify in this paper strongly suggests that 
 the Virasoro constraint conjecture of \cite{art:Eguchi1, art:Eguchi2}
 is a reflection of the combinatorial structure of building 
 the domain Riemann 
 surface from simpler objects such as pairs of pants or three punctured
 spheres.

 Although it is more than 15 years old, the Witten-Kontsevich
 theory \cite{art:Witten2dGravity, art:Kontsevich}
  has never lost its place as one of the most beautiful and
 prime theories in the study of algebraic curves and their moduli
 spaces. The theory provides a complete computational method
 for  all intersection numbers of the tautological 
 cotangent classes (the $\psi$-classes) defined on the moduli 
 space $\CompactModuli_{g,n}$ of stable algebraic curves 
 of genus $g$ with $n$ marked points.
    Recently several new proofs
have appeared \cite{art:OkounkovPandharipande, art:MirzakhaniIT, art:KimLiu, art:KazarianLando}. We note that  all these new proofs are based on very different 
 ideas and techniques, including random matrix theory, random
 graphs, Hurwitz theory, representation theory of 
 symmetric groups, symplectic geometry and
 hyperbolic geometry. 

 The mystery of the Witten-Kontsevich theory has been the 
 following question:
 where does the KdV equation, and also the Virasoro
 constraint condition, come from? Once we accept  the 
 Kontsevich matrix model expression of the generating 
 function of all cotangent class intersections, then both 
 the KdV and the Virasoro are an easy corollary of the
 analysis of matrix integrals. Thus the real question is:
 where do these structures appear in the geometry of
 moduli spaces of algebraic curves?
 
 The insight of some of the
 the new proofs \cite{art:KimLiu, art:KazarianLando}, which do not
 rely on matrix integrals but rather use the counting of ramified
 overings of $\mathbb{P}^1$,
 is that  the
 KdV equation is a direct consequence of the
 \emph{cut and join} mechanism of \cite{art:GouldenJacksonV}.

 The proof \cite{art:MirzakhaniIT} due to Mirzakhani
 utilizes hyperbolic geometry and has a markedly different nature from
the  others, whose origins are  rooted in algebraic geometry. 
Mirzakhani's work concerns  the Weil-Petersson volume of the moduli space of \emph{bordered}
 Riemann surfaces $\Moduli_{g,n}(\MathVector{L})$, where $\MathVector{L} = (L_1, \ldots, L_n)$ specifies the geodesic lengths of the boundaries of Riemann surfaces. 
 Here the moduli space is equipped with the
 structure of  a differentiable orbifold
 realized as the quotient of the Teichm\"uller space by the action 
 of a mapping class group. 
 Mirzakhani shows that these volumes satisfy a 
 recursion relation, and that
 in the limit $\MathVector{L} \rightarrow \infty$ her recursion formula recovers the Virasoro constraint condition for the generating function of $\psi$-class intersection numbers of $\CompactModuli_{g,n}$, or
 equivalently, the generating function of the Gromov-Witten invariants
 of a point. 
 A striking theorem of \cite{art:MirzakhaniIT} relates, via
 the method of symplectic reduction,  the 
 Weil-Petersson volume of  $\Moduli_{g,n}(\MathVector{L})$ 
 and the intersection numbers involving both the first Mumford
 class $\kappa_1$ and the $\psi$-classes on 
$\CompactModuli_{g,n}$. As a 
 consequence, she proves that the volume 
 $\Vol\big(\Moduli_{g,n}(\MathVector{L})\big)$, after an 
 appropriate normalization with powers of $\pi$, is a 
 polynomial in $\MathVector{L}$ with rational coefficients.  
 
 Since there is no particular reason to believe that there
 should be a direct relation between the Weil-Petersson volume of the
 moduli spaces of bordered Riemann surfaces and Virasoro
 constraint condition, our discovery suggests the existence of 
 another, more algebraic, point of view in the Mirzakhani theory.

   Our Virasoro structure also bears an interesting consequence: 
 it leads to the natural normalization of the
 Weil-Petersson volume of the moduli spaces of 
 bordered (or unbordered)  Riemann surfaces.  
 Although the geometric orbifold picture and the algebraic stack
 picture give the same moduli space for most of the cases, there is one
 exception:
 the moduli space of one-pointed 
 stable elliptic curves $\CompactModuli_{1,1}$. If we define this space as 
 an orbifold, then its canonical Weil-Petersson
 volume is $\zeta(2) = \pi^2/6$. On the other hand, the Virasoro
 constraint condition dictates that we need to have
 $$
 \Vol(\CompactModuli_{1,1}) = \frac{\zeta(2)}{2}
 $$
 as its canonical symplectic volume. This makes sense
 if we consider $\CompactModuli_{1,1}$ as an algebraic stack. 
 The factor $2$ difference is due to the fact that every 
 elliptic curve with one marked point possesses a 
 $\mathbb{Z}/2\mathbb{Z}$ automorphism. It is remarkable 
 that even a purely hyperbolic geometry argument leads us to 
 this stack picture.

To summarize our main results, let us consider the rational volume of 
$\Moduli_{g,n}(\MathVector{L})$ defined by
\begin{equation*}
	v_{g,n}(\MathVector{L}) \stackrel{\text{def}}{=} 
		\frac{\Vol\big(\Moduli_{g,n}(\MathVector{L})\big)}{2^d\pi^{2d} } ,
\end{equation*}
where $d= 3g-3+n$ and
\begin{equation*}
	\Vol\big(\Moduli_{g,n}(\MathVector{L})\big) \stackrel{\text{def} }{= } \int_{\Moduli_{g,n}(\MathVector{L})} \frac{\omega_{WP}^d }{d! }
\end{equation*}
is the Weil-Petersson volume of $\Moduli_{g,n}(\MathVector{L})$.
Then the Mirzakhani recursion formula reads 
\begin{align*}
	v_{g,n}(\MathVector{L}) &= \frac{2}{L_1} \int_{0}^{L_1} \int_{0}^{\infty} \int_{0}^{\infty}
		xy K(x+y, t) v_{g-1,n+1}(x,y,\MathVector{L}_{\hat{1}}) dxdydt \\
	&+ \frac{2}{L_1} \sum_{\substack{g_1+g_2=g\\ \mathcal{I} \coprod \mathcal{J} = \{2,\ldots,n\}} }
		\int_{0}^{L_1}\int_{0}^{\infty}\int_{0}^{\infty} xy K(x+y, t) 
		v_{g_1,n_1}(x, \MathVector{L}_{\mathcal{I}})  \\[-2mm]
		& \hspace{53mm} \times v_{g_2,n_2}
		(y, \MathVector{L}_{\mathcal{J}})
		dxdydt \\[2mm]
	&+ \frac{1}{L_1} \sum_{j=2}^{n} \int_{0}^{L_1} \int_{0}^{\infty} x 
		\left( K(x, t+L_j) + K(x, t-L_j)\right) \\
		&\hspace{30mm} \times v_{g, n-1}(x, \MathVector{L}_{\{\widehat{1,j}\} }) dxdt,
\end{align*}
where the kernel function of the integral transform is given by
\begin{equation*}
	K(x,t) = \frac{1}{1 + e^{\pi(x+t)}} + \frac{1}{1 + e^{\pi(x-t)}},
\end{equation*}
and the symbol $\widehat{\hspace{10pt}}$ indicates the complement of the indices.
Recall  that our normalized Weil-Petersson volume is a polynomial in $\MathVector{L}$ with coefficients given by intersection numbers of $\kappa_1$ and $\psi$-classes:
\begin{equation*}
	v_{g,n}(\MathVector{L}) 
			= \sum_{\substack{d_0 + \cdots + d_n\\ = d}} \prod_{i=0}^{n}\frac{1}{d_i! } 
		\intersect{\kappa_1^{d_0} \prod\tau_{d_i} }_{g,n} \prod_{i=1}^{\infty}
		L_{i}^{2d_i}.
\end{equation*}
Instead of defining our generating function directly from these rational 
volumes, let us consider the generating function of the
mixed $\kappa_1$ and $\psi$-class intersections  
\begin{align*}
	G(s, t_0, t_1, t_2, \ldots) &\stackrel{\text{def}}{=}
		\sum_{g} \intersect{e^{s\kappa_1 + \sum t_i \tau_i} }_{g} \\
		&= \sum_{g} \sum_{m, \{n_i\}} \intersect{\kappa_1^m \tau_0^{n_0}\tau_1^{n_1}
			\cdots }_{g} 
		\frac{s^m}{m!} \prod_{i=0}^{\infty} \frac{t_{i}^{n_i}}{n_i! }.
\end{align*}
 The main results of the present paper are the following
 \emph{differential} version of the integral recursion formula.
\begin{theorem} For every $k\ge -1$, let us define
 \begin{multline*}
 	\Virasoro_k = -\frac{1}{2} \sum_{i=0}^{\infty} (2(i+k)+3)!! \alpha_{i} s^i 
		\frac{\partial }{\partial t_{i+k+1} } 
		+ \frac{1}{2} \sum_{j=0}^{\infty} \frac{(2(j+k)+1)!! }{(2j-1)!! } t_j
		\frac{\partial }{\partial t_{j+k} } \\
	+ \frac{1}{4} \sum_{\substack{d_1 + d_2 = k-1 \\ d_1, d_2 \geq 0 } }
		(2d_1 + 1)!! (2d_2 + 1)!! \frac{\partial^2 }{\partial t_{d_1} \partial t_{d_2}}
		+ \frac{\delta_{k,-1}t_0^2}{4} + \frac{\delta_{k,0} }{48},
 \end{multline*}
 where $\alpha_{i} = \frac{(-2)^i }{(2i+1)! }$. Then we have:
\begin{enumerate}
\item
	The operators $V_k$ satisfy Virasoro relations
\begin{equation*}
	[V_{n}, V_{m}] = (n-m)V_{n+m}.
\end{equation*}
\item
The function $\exp(G)$ satisfies the Virasoro constraint condition
\begin{equation*}
	V_k \exp(G) = 0 \quad \text{for $k\ge -1$}.
\end{equation*}
\end{enumerate}
Moreover, these properties uniquely determine $G$ and enable one to calculate all coefficients of the expansion.  
Since $G$ contains all
information of the rational volumes $v_{g,n}(\MathVector{L})$, we conclude that the
Virasoro constraint condition is indeed equivalent to the Mirzakhani
recursion relation.
\label{thm:IntroVirasoro}
\end{theorem}

Since the generating function for $\psi$-class intersections
\begin{align*}
	F(t_0, t_1, \ldots) &= \sum_{g} \intersect{e^{\sum \tau_i t_i} }_{g} \\
		&= \sum_{g} \sum_{\{n_{*} \} } \intersect{\prod \tau_{i}^{n_i} }_{g} \prod \frac{t_i^{n_i} }{n_i! }
\end{align*}
is a solution of the KdV hierarchy, it is natural to ask if $G$ satisfies any integrable equations. Indeed, we prove the following.
\begin{theorem} The function
$\exp(G)$ is a $\tau$-function for the KdV hierarchy for any fixed value of $s$. In fact, we have an explicit relation
 	\begin{equation}
		G(s, t_0, t_1, \ldots) = F(t_0, t_1, t_2 + \gamma_2, t_3 + \gamma_3, \ldots),
		\label{eqn:GeneratingEquiv}
	\end{equation}
where $\gamma_i = \frac{(-1)^i}{(2i+1)i!}s^{i-1}$. 
\label{thm:IntroKdV}
\end{theorem}
\noindent We remark that it is well known to algebraic geometers that generating functions $F$ and $G$ contain the same information \cite{art:ArbarelloCornalba2, art:ArbarelloCornalba1,  art:Faber, art:KMZ, art:ManinZograf, art:Zograf}.

An important consequence of Theorem~\ref{thm:IntroKdV} is that $G$ is also completely determined by the property of being a $\tau$-function, together with the string equation $V_{-1}\exp(G)=0$.  It is fruitful to think of the string equation as being the initial condition for the KdV flow. Since $G$ is determined, we note that Theorem~\ref{thm:IntroKdV} is again equivalent to Mirzakhani's recursion formula.

Here we recall that in the theory of integrable systems, every 
variable has a weighted degree so that all natural operators have
homogenous weights. Coming from the KdV equations, we 
assign $\deg t_j = 2j+1$. The quantity $\gamma_j$ has the same 
degree, which defines that $\deg s^j = 2j + 3$. The Virasoro
operator $V_k$ then has homogenous degree $-2k$ for every
$k\ge -1$. 
Another way to view the degree of $s^i$ comes from the generalized Kontsevich integral
\begin{equation*}
	\log \int_{\mathcal{H}_N } e^{i\sum_{j=0}^{\infty} (-\frac{1}{2} )^j s_j \frac{\trace X^{2j+1} }{2j+1 }}
	e^{-\frac{\trace (X^2\Lambda) }{2 } }dX.
\end{equation*}
As indicated in the work of Mondello \cite{art:Mondello}, there should be a substitution $s_j = c_j s^{j-1 }$ which transforms the asymptotic expansion of the integral into the generating function $G$. Since $s_j$ has degree $2j+1$, we confirm that $s^j$ must have degree $2j+3$.
As well, we should point out that it is quite natural for (\ref{eqn:GeneratingEquiv}) to leave variables $t_0$ and $t_1$ unchanged. The reason is that in any expression relating intersections involving $\kappa$ classes to those involving $\tau$ terms alone, $\tau_0$ and $\tau_1$ never make an appearance.

A few of the natural questions that crop up from this work are: 
\begin{enumerate}
\item Is there a matrix integral expression for the function $G$? 
\item Is it possible to prove that $G$ is a solution to the KdV hierarchy without appealing to the Witten-Kontsevich theorem? 
\item What is the direct geometric connection between the cut and join
mechanism and the Mirzakhani recursion?
\end{enumerate}
We note that the essence of the original Virasoro constraint conjecture
is that the generating function of Gromov-Witten invariants \emph{should}
have a matrix integral expression. The analysis of matrix integrals
\cite{art:AdlervanMoerbecke} indicates that once a matrix integral
formula is established, the Virasoro constraints and integrable systems 
of KdV type are obvious consequences. Since the ribbon graph
expansion method provides a powerful tool to matrix integrals, the very
existence of both the KdV equations and the Virasoro constraint for the
Weil-Petersson volume of the moduli spaces of bordered Riemann
surfaces points to a matrix model expression and ribbon 
graph interpretation of the Mirzakhani
formulas. 
These questions, however,  are beyond the scope of our present work.

This paper is organized as follows.
In section~\ref{sect:Mirzakhani} we review the
work of Mirzakhani  \cite{art:MirzakhaniWP, art:MirzakhaniIT}.  
Since the Virasoro structure very delicately depends on all the subtle points 
of the theory, we provide a detailed 
discussion on some of the key ingredients of the
work, including the case of genus one with one boundary, precise 
combinatorial description of cutting a surface along geodesics, and the
choice of a canonical orientation of the circle bundle when the
Duistermaat-Heckman formula is applied to the extended moduli
spaces. 
Section~\ref{sect:Virasoro} gives a proof of Theorem~\ref{thm:IntroVirasoro}. Finally, in section~\ref{sect:IntegrableSystem} we prove Theorem~\ref{thm:IntroKdV}.

The second author would like to thank Greg Kuperberg and  Albert Schwarz for helpful conversations. 

 \section{Mirzakhani's Recursion Relation}
\label{sect:Mirzakhani}

\subsection{Notations}

Since the orbifold picture and the stack structure are
the same for the moduli spaces of algebraic curves except
for genus $1$ with one marked point, we employ the orbifold
view point throughout the paper. As mentioned above, however, 
when we interpret the canonical volume of the moduli spaces,
we need to use the stack picture. 

Let $\Moduli_{g,n}$ denote the moduli space of smooth algebraic
curves, or equivalently, the moduli orbifold 
consisting of finite area hyperbolic metrics on a surface,
 of topological type $(g,n)$.
 A surface of type $(g,n)$ is a surface with genus $g$ and $n$ punctures. 
 Since we  are interested in the stable, noncompact case,  we
 impose $2g-2+n > 0$ and $n>0$ throughout this paper. When referring to the underlying topological type of the surface we will consistently employ the notation $S_{g,n}$.  We will also use the notation $\Moduli_{S}$ for the moduli space of surfaces of topological type $S$. 
  
  Mirzakhani's breathtaking theory is about the moduli space
 $\Moduli_{g,n}(\MathVector{L})$ of genus $g$ hyperbolic surfaces with $n$ geodesic boundary components of specified length $\MathVector{L} = (L_1, \ldots, L_n)$. This space relates to the algebro-geometric 
 moduli space via the equality $\Moduli_{g,n} = \Moduli_{g,n}(0)$.  
 The moduli space of bordered Riemann surfaces is defined
 as an orbifold 
 \begin{equation*}
 	\Moduli_{g,n}(\MathVector{L}) = \Teichmuller_{g,n}(\MathVector{L}) / \MCG_{g,n},
 \end{equation*} 
 where $\MCG_{g,n}$ is the mapping class group of the surface of type $(g,n)$, i.e., the set of isotopy classes of diffeomorphisms which preserve the boundaries setwise, and $\Teichmuller_{g,n}(\MathVector{L})$ is
 the Teichm\"uller space. The Deligne-Munford type compactification 
 of this moduli space is obtained by pinching non-trivial cycles.

 The tautological classes we consider in this
 paper are the $\kappa_1$ and $\psi$ classes. Let
  $$
 \pi: \CompactModuli_{g,n+1}
 (= \mathcal{C}_{g,n}) \longrightarrow \CompactModuli_{g,n}
 $$
 be the forgetful morphism
  which forgets the $n+1$-st marked point, and
  $$
  \sigma_{i}(C, x_1, \ldots, x_n) = x_i \in C, \qquad i=1,2,\dots,n
  $$
  its canonical sections. 
 	We denote by $\omega_{\mathcal{C} / \Moduli}$ the relative dualizing sheaf,  and let $\mathcal{D}_i = \sigma_i(\CompactModuli_{g,n})$, which
	is a divisor
	in $\CompactModuli_{g,n+1}$. The tautological classes are defined
	by
\begin{align*}
	\mathcal{L}_i &= \sigma_i^*(\omega_{\mathcal{C}/\Moduli }), \\
	\psi_i &= c_1(\mathcal{L}_i), \\
	\kappa_1 &= \pi_* \biggl( c_1\bigl(\omega_{\mathcal{C}/\Moduli}({\textstyle \sum} \mathcal{O}(\mathcal{D}_i) \bigr)^{2}  \biggr).
\end{align*}
We are interested in the intersection numbers 
\begin{equation*}
	\intersect{\kappa_1^m \tau_{d_1} \cdots \tau_{d_n} }_{g}
		= \int_{\CompactModuli_{g,n}} \kappa_{1}^m 
			\psi_{1}^{d_1} \cdots \psi_{n}^{d_n}.
\end{equation*}

The class $\kappa_1$ has a nice geometric interpretation, coming from the symplectic structure of $\Moduli_{g,n}$. The Fenchel-Nielsen coordinates
are associated with a pair of pants decompostion of the surface $S_{g,n}$, which is a disjoint set of simple closed curves $\Gamma = \{\gamma_1, \ldots, \gamma_d\}$ such that $S_{g,n} \setminus \Gamma$ is a disjoint union of pairs of pants (triply punctured spheres). Since pairs of pants have no moduli (they are uniquely fixed after specifying the boundary lengths), all that remains to recover the original hyperbolic structure is to specify how the matching geodesics are glued together. Hence for every curve $\gamma_i$ in the pair of pants decomposition, one has the freedom of two parameters $(l_i, \tau_i)$ where $l_i$ is the length of the curve and $\tau_i$ is the twist parameter. These coordinates give an isomorphism $\Teichmuller_{g,n} = \Reals_{+}^d \times \Reals^d$. At the moduli space level we must quotient out by the mapping class group action. Since one full twist around a curve is a Dehn twist (element of $\MCG_{g,n}$), we get local coordinates of the form $\Reals_{+}^d \times (S^1)^d$. 

In Fenchel-Nielsen coordinates, the Weil-Petersson form is in Darboux coordinates \cite{art:WolpertFNtwist}:
\begin{equation*}
	\omega_{WP} = \sum dl_i \wedge d\tau_i.
\end{equation*}
This is a closed, nondegenerate 2-form on $\Teichmuller_{g,n}$ which is invariant under the action of the mapping class group, hence gives a well defined symplectic form on $\Moduli_{g,n}$. Note that by Wolpert~\cite{art:WolpertWP} the Weil-Petersson form extends as a closed current on $\CompactModuli_{g,n}$.  In particular, the Weil-Petersson volume
\begin{equation*}
	\Vol_{g,n}(\MathVector{L}) \stackrel{\text{def} }{= } \int_{\Moduli_{g,n}(\MathVector{L})} \frac{\omega_{WP}^d }{d! }
\end{equation*}
is a finite quantity. (Note that Wolpert defines the Weil-Petersson form to be half of the above expression. Our convention is adopted from the algebraic geometry community \cite{art:ArbarelloCornalba2, art:KMZ, art:Zograf}.)
The relation to tautological classes is provided by the well-known 
formula	
$$\omega_{WP} = 2\pi^2 \kappa_1.
$$

 For use in the sequel, we note that the vector field generated by a Fenchel-Nielsen twist about a simple closed geodesic is symplectically dual to the length of the geodesic. A Fenchel-Nielsen twist is defined by cutting the surface along the curve, twisting one component with respect to the other and than regluing. As a formula, we have
 \begin{equation*}
 	\omega_{WP}(\, \cdot\,  , \frac{\partial}{\partial \tau_i}) = dl_i.
 \end{equation*}
 

\subsection{McShane's identity}

A crucial step in Mirzakhani's program \cite{art:MirzakhaniWP, art:MirzakhaniIT} is to use McShane's identity to write a constant
function on the moduli space
 as a sum over mapping class group orbits of simple closed curves. 
To state Mirzakhani's generalization of McShane's identity, we introduce the following notation for an arbitrary hyperbolic surface $X$ with boundaries $(\beta_1, \ldots, \beta_n)$ of length $(L_1, \ldots, L_n)$:
	$\mathcal{I}_j$ denotes the set of 
	simple closed geodesics $\gamma$ such
	that $(\beta_1,\beta_j, \gamma)$ bound a pair of pants; and
	$\mathcal{J}$ the set of  pairs of simple closed geodescis $(\alpha_1,\alpha_2)$ such that $(\beta_1, \alpha_1, \alpha_2)$ bound a pair of pants.
Using the functions
\begin{align*}
	\mathcal{D}(x,y,z) &= 2\log \left(
		\frac{e^{x/2} + e^{(y+z)/2} }{e^{-x/2} + e^{(y+z)/2 } } \right) \text{ and}\\
	\mathcal{R}(x,y,z) &= x - \log \left(
		\frac{\cosh\tfrac{y}{2} + \cosh\tfrac{x+z }{2 } }{
		   \cosh\tfrac{y }{2} + \cosh\tfrac{x-z }{2} } \right),
\end{align*}
 Mirzakhani proves
\begin{theorem}[Mirzakhani \cite{art:MirzakhaniWP}]
	For $X$ as above, we have 
	\begin{equation*}
	   L_1 = \sum_{(\alpha_1,\alpha_2)\in\mathcal{J} } \mathcal{D}\bigl(L_1, l(\alpha_1),
	   	l(\alpha_2)\bigr)
	+ \sum_{j=2}^{n} \sum_{\gamma\in\mathcal{I}_j } \mathcal{R} \bigl(L_1,
		L_j, l(\gamma) \bigr).
	\end{equation*}
\end{theorem}

\subsection{Integration over the moduli spaces}

The idea is to find a fundamental domain for a particular cover of $\Moduli_{g,n}$, enabling one to integrate functions over this covering space. Then a specific class of functions defined on the moduli space (such as those arising from McShane's identity) can be lifted to this cover and integrated. To that end, let $\Gamma = \{\gamma_1, \ldots , \gamma_n\}$ be a collection of disjoint simple closed curves on the surface $S_{g,n}$, where $S_{g,n}$ is the underlying topology of the hyperbolic surface $X\in\Moduli_{g,n}$. We define
\begin{equation*}
	\Stab \Gamma = \cap \Stab \gamma_i = \{f\in\MCG_{g,n}\, | \, f(\gamma_i) = 	
		\gamma_i, \, \text{for $i=1, \ldots, n$} \},
\end{equation*}
and set
\begin{align*}
	\Moduli_{g,n}^{\Gamma} &= \Teichmuller_{g,n} / \Stab\Gamma \\
	 &= \{ (X, \eta_1, \ldots, \eta_n)\, |\, X\in\Moduli_{g,n},\ 
	 \text{$\eta_i$ is a simple closed geodesic in $\MCG\gamma_i$ }  \}.
\end{align*}

Note that as a quotient of Teichm\"uller space, $\Moduli_{g,n}^{\Gamma}$ inherits the Weil-Petersson symplectic form. Hence we can talk about integration over $\Moduli_{g,n}^{\Gamma}$ with respect to the symplectic volume form. The advantage of integration on $\Moduli^{\Gamma}_{g,n}$ as opposed to the usual moduli space is that we can exploit the existence of a hamiltonian torus action. In fact, by a result of Wolpert, the vector field generated by a Fenchel-Nielsen twist along a geodesic is symplectically dual to the length function of the geodesic (as a function on Teichm\"uller space). The space $\Moduli_{g,n}^{\Gamma}$ is the intermediate covering space on which the circle actions on $\{\gamma_1, \cdots, \gamma_n\}$ descend. The problem with attempting to construct such a circle action on $\Moduli_{g,n}$ is that there is no well defined notion of a geodesic curve on an element of moduli space; the best one can obtain is a mapping class group orbit of curves. 

\begin{figure}
	\begin{pspicture}(7,2.9)
	\rput[bl](0.3,0){\epsfbox{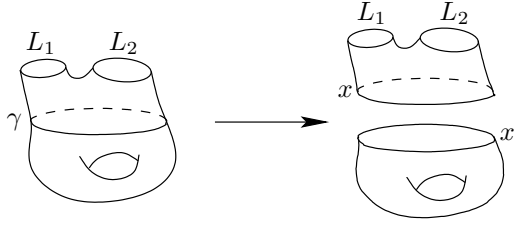}}
	\uput{0.05}[l](0.4,1.3){$\gamma$}
	\uput{0.05}[u](0.6,2.15){$L_1$}	\uput{0.05}[u](1.7, 2.15){$L_2$}
	\uput{0.05}[l](4.8,1.7){$x$} 		\uput{0.05}[r](6.65, 1.13){$x$}
	\uput{0.05}[u](5, 2.55){$L_1$} 		\uput{0.05}[u](6.1, 2.55){$L_2$}
	\end{pspicture}
	\caption{Decomposing a surface}
	\label{fig:SurfaceDecomp}
\end{figure}

Thus we have the moment map for the torus action:
\begin{align*}
	\MathVector{l}: \Moduli_{g,n}^{\Gamma} &\rightarrow \Reals_{+}^{n} \\
	(X, \MathVector{\eta}) &\mapsto \bigl(l(\eta_1), \ldots, l(\eta_n)\bigr).
\end{align*}
Hence we see that $\MathVector{l}^{-1}(\MathVector{x}) / T^n$ is symplectomorphic to $\Moduli_{S_{g,n} \setminus \Gamma}(\MathVector{L}, \MathVector{x}, \MathVector{x})$. Recall that $\Moduli_{S_{g,n} \setminus \Gamma}$ is the moduli space with underlying topological type the (possibly disconnected) surface $S_{g,n} \setminus \Gamma$ and with boundary lengths following the rule outlined in Figure~\ref{fig:SurfaceDecomp}.

The most straightforward way to prove the above assertion is to take a pair of pants decomposition for the surface $S_{g,n}$ which contains the curves $\Gamma$.  Then $\MathVector{l}^{-1}(\MathVector{x})$ fixes the lengths of the geodesics $\eta_1, \ldots, \eta_n$, while quotienting by the torus action removes the twist variable from these curves. This gives the diffeomorphism. The symplectic equivalence follows immediately from the Fenchel-Nielsen coordinate expression for the Weil-Petersson form.
What emerges is an exceptionally clear local picture for the space $\Moduli_{g,n}^{\Gamma}$. In fact, it is a fibre bundle over $\Reals_{+}^{n}$ where the fibres are (locally) equal to the product of a torus  and $\Moduli_{S_{g,n} \setminus \Gamma}(\MathVector{L}, \MathVector{x}, \MathVector{x})$. 

Consider a map
\begin{equation*}
	f:\Moduli_{g,n}^{\Gamma} \rightarrow \Reals,
\end{equation*}
which is a function of the lengths of the marked geodesics $l(\eta_i)$. In other words, $f(X, \eta_1, \ldots, \eta_n) = f\bigl(l(\eta_1), \ldots, l(\eta_n)\bigr)$. By the previously discussed decomposition of $\Moduli_{g,n}^{\Gamma}$ we can write
\begin{equation}
	\int_{\Moduli_{g,n}^{\Gamma} } f\bigl(\MathVector{l}(\MathVector{\eta})\bigr)
	e^{\omega_{WP}(\MathVector{L})} = 
	\int_{\Reals_{+}^{n} } f(\MathVector{x}) \Vol_{S_{g,n} \setminus \Gamma} 
	(\MathVector{L}, \MathVector{x}, \MathVector{x}) 
	\MathVector{x}\cdot d\MathVector{x}.
\end{equation}
Here $e^{\omega_{WP}}$ means we are integrating over the maximal power of the Weil-Petersson form $\frac{\omega_{WP}^{d} }{d! }$,  $d = 3g-3+n$.
Note that if $S_{g,n} \setminus\Gamma$ is disconnected then $\Moduli_{S_{g,n} \setminus\Gamma}$ is the direct product of the component moduli spaces, with the volume being the product of each.

\subsection{Volume calculation}
To relate the above discussion to integration on the moduli space, 
Mirzakhani uses her generalized McShane identity. As a lead in to the main results, consider the following simplified situation. Suppose $\gamma$ is a simple closed curve on $S_{g,n}$, with $\MCG\gamma$ its mapping class group orbit. Then given any hyperbolic structure $X$ on $S_{g,n}$, every $\alpha\in\MCG\gamma$ has a unique geodesic in its isotopy class. Denote $l_{X}(\alpha)$ the corresponding geodesic length. Hence for any function $f:\Reals_+ \rightarrow \Reals$ and  thinking of $X\in[X]$ as a representative of an element of $\Moduli_{g,n}$ we have the following well defined function on $\Moduli_{g,n}$:
\begin{align*}
	f^{\gamma}: \Moduli_{g,n} &\rightarrow \Reals \\
	[X] &\mapsto \sum_{\alpha\in\MCG\gamma} l_{X}(\alpha).
\end{align*}
One can easily check that this function does not depend on the choice of representative $X\in[X]$. However, it is not a priori clear whether or not $f^{\gamma}$ will be a convergent sum. At minimum one requires $\lim_{x\rightarrow\infty}f(x) = 0$.
We can similarly define a function
\begin{align*}
	\tilde{f}^{\gamma}: \Moduli_{g,n}^{\gamma} &\rightarrow \Reals \\
\intertext{by the rule}
		\tilde{f}^{\gamma}(X, \eta) &= f(l_X(\eta)),
\end{align*}
which gives the relation
\begin{equation*}
	f^{\gamma}(X) = \sum_{(X,\eta)\in\pi^{-1}(X)} \tilde{f}^{\gamma}(X,\eta).
\end{equation*}
In particular, since the pullback of the Weil-Petersson form is the Weil-Petersson form on the cover, we have
\begin{equation*}
	\int_{\Moduli_{g,n}} f^{\gamma} e^{\omega_{WP}} 
		= \int_{\Moduli_{g,n}^{\gamma}} \tilde{f}^{\gamma} e^{\omega_{WP}}.
\end{equation*}

\begin{figure}
	\begin{pspicture}(9.7, 4.2)
	\rput[bl](0,0){\epsfbox{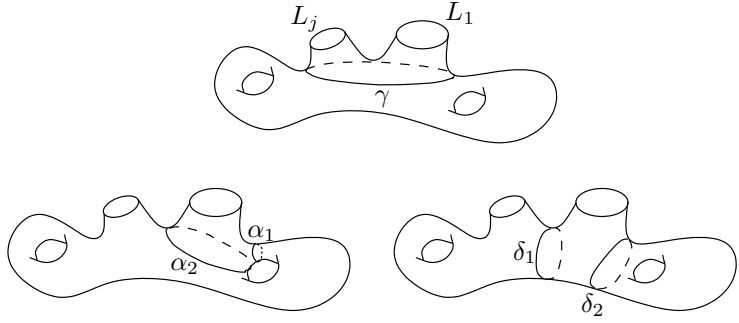}}
	\uput{0.05}[dl](2.6,0.75){$\alpha_2$} \uput{0.05}[u](3.4,0.9){$\alpha_1$}
	\uput{0.05}[ul](4.2,3.7){$L_j$} \uput{0.05}[ur](5.8,3.8){$L_1$}
	\uput{0.05}[d](5, 3){$\gamma$}
	\uput{0.05}[l](7.1, 0.8){$\delta_1$}  \uput{0.05}[d](7.8,0.3){$\delta_2$}
	\end{pspicture}
	\caption{Removing a pair of pants from a surface}
	\label{fig:PantsRemoval}
\end{figure}

Note that for any curve $\gamma\in\mathcal{I}_j$, we have $\mathcal{I}_j = \MCG\gamma$. This follows because two curves are in the same orbit of the mapping class group if and only if the surfaces obtained by cutting along the curves are homeomorphic, with a homeomorphism preserving the boundary components setwise. The homeomorphism will extend continuously to the curves to give the map of the entire surface. This tells us that the set $\mathcal{J}$ is not the orbit of a single pair of curves $(\alpha_1,\alpha_2)\in\mathcal{J}$. In fact, we further refine this set of curves as follows. For any $(\alpha_1,\alpha_2)\in\mathcal{J}$ set $P(\beta_1, \alpha_1, \alpha_2) \subset S_{g,n}$ to be the pair of pants bounded by the curves $\beta_1, \alpha_1, \alpha_2$. Now we define (see Figure~\ref{fig:PantsRemoval})
\begin{align*}
	\mathcal{J}_{\text{conn}} &= \{ (\alpha_1, \alpha_2)\in\mathcal{J}\, |\,
		\text{$S_{g,n} \setminus P(\beta_1, \alpha_1,\alpha_2)$ is connected }   \}\\
	\mathcal{J}_{g_1, \{i_1,\ldots , i_{n_1} \}} &=  \{ (\alpha_1, \alpha_2)\in
		\mathcal{J}\, |\, \text{$S_{g,n} \setminus P(\beta_1, \alpha_1,\alpha_2)$
		breaks into 2 pieces } \\
	& \quad\quad \text{one of which is a surface of type $(g_1,n_1+1)$ with  } \\
	& \quad\quad \text{boundary $(\alpha_i, \beta_{i_1}, \ldots, \beta_{i_{n_1}}) $ } \}.
\end{align*}
Other than the obvious identification 
\begin{equation*}
	\mathcal{J}_{g_1, \{i_1,\ldots , i_{n_1} \}} =
 	\mathcal{J}_{g-g_1, \{1,\ldots, n\} \setminus \{i_1,\ldots , i_{n_1} \}},
\end{equation*}
these subsets form disjoint orbits under the mapping class group. Moreover $\MCG_{g,n}$ acts transitively on each set.

Hence we can write the Mirzakhani-McShane identity in the following form
\begin{align*}
	L_1 &= \frac{1}{2} \sum_{\substack{g_1+g_2 = g \\ \mathcal{A} \coprod 
		\mathcal{B} = \\ \{2,\ldots,n\}} } 
		\sum_{(\alpha_1,\alpha_2)\in \mathcal{J}_{g_1, \mathcal{A} }} 
		\mathcal{D}\bigl(L_1, l(\alpha_1 ), l(\alpha_2 )\bigr) \\
		&+ \sum_{(\delta_1,\delta_2)\in \mathcal{J}_{\text{conn}} } 
		  \mathcal{D}\bigl(L_1, l(\delta_1 ), l(\delta_2 )\bigr) \\
		 &+ \sum_{j=2 }^{n} \sum_{\gamma\in\mathcal{I}_j } 
		 \mathcal{R} \bigl(L_1,L_j, l(\gamma) \bigr).
\end{align*}
There is a slight inaccuracy - we undercount by half for terms with  $n=1$ and $g_1 = g_2$. However, we will see in a moment that this makes further calculations somewhat simpler.

Each of the terms in the above sum can be lifted to a function on an appropriate cover $\Moduli_{g,n}^{\Gamma}$. We see that
\begin{align*}
	\int_{\Moduli_{g,n}(\MathVector{L})} L_{1} e^{\omega_{WP} } 
	&= 
		\frac{1}{2} \sum_{\substack{g_1+g_2 = g \\ \mathcal{A} \coprod \mathcal{B}
		 = \\ \{2, \ldots, n\} } }  \sum_{(\alpha_1,\alpha_2) \in 
		 	\mathcal{J}_{g_1, \mathcal{A} } }
		\int_{\Moduli_{g,n}(\MathVector{L} ) } \mathcal{D} 
		(L_1, l(\alpha_1), l(\alpha_2) )e^{\omega_{WP} }
	 \\
	 &+ \sum_{(\delta_1,\delta_2 ) \in \mathcal{J}_{\text{conn}} } \int_{\Moduli_{g,n}	(\MathVector{L} ) } \mathcal{D} (L_1, l(\delta_1), l(\delta_2)) e^{\omega_{WP}} \\
	 &+ \sum_{j=2}^{n} \sum_{\gamma\in\mathcal{I}_j } 
	 \int_{\Moduli_{g,n}(\MathVector{L} ) } \mathcal{R} 
	 (L_1, L_j, l(\gamma) )e^{\omega_{WP}}.
\end{align*}
	Hence
\begin{align*}
	L_1 \Vol_{g,n}(\MathVector{L} )
	&=  \frac{1}{2} \sum_{\substack{g_1+g_2 = g \\ \mathcal{A} \coprod \mathcal{B}
		 = \{2, \ldots, n\} } }  \int_{\Moduli_{g,n}^{\{\alpha_1,\alpha_2 \} } }
		\mathcal{D}(L_1,l(\eta_1),l(\eta_2)) e^{\omega_{WP} } 
	 \\
	 &+ \frac{1}{2} \int_{\Moduli_{g,n}^{\{\delta_1, \delta_2\} } } 
	 	\mathcal{D} (L_1, l(\eta_1), l(\eta_2)) e^{\omega_{WP} } \\
	&+ \sum_{j=2}^{n} \int_{\Moduli_{g,n}^{\gamma } } \mathcal{R}
		(L_1, L_j, l(\gamma) ) e^{\omega_{WP} }.
\end{align*}

Note that the factor $\frac{1}{2}$ that appears in front of the second term in the sum is needed because the mapping class group orbit double counts the set of curves $(\delta_1, \delta_2 )$. In other words, there is a diffeomorphism exchanging $\delta_1$ and $\delta_2$. Similarly, the factor $\frac{1}{2}$ discrepancy for the term in the first sum with $g_1 = g_2$ and $n=1$ has now disappeared and the presented sum is unambiguously correct. Applying the results of the previous section we see that 
\begin{align*}
	L_1 \Vol_{g,n}(\MathVector{L} )
	&= \frac{1}{2} \sum_{\substack{g_1 + g_2 = g \\ \mathcal{A} \coprod \mathcal{B} } } 
		\int_{\Reals_{+}^2} x y \mathcal{D} (L_1, x ,y)
		\Vol_{g_1, n_1 }(x,\MathVector{L}_{\mathcal{A}} )
		\Vol_{g_2,n_2} (y, \MathVector{L}_{\mathcal{B}} ) dxdy
	\\
	&+ \frac{1}{2} \int_{\Reals_{+}^{2} } xy \mathcal{D}(L_1,x,y ) 
		\Vol_{g-1,n+1}(x,y,\MathVector{L }_{\hat{1} } ) dxdy \\
	&+ \sum_{j=2}^{n} \int_{\Reals_+ } x \mathcal{R} (L_1,L_j, x)
		\Vol_{g,n-1 }(x, \MathVector{L}_{\widehat{1,j} } )dx.
\end{align*}
For any subset $\mathcal{A} = \{i_1, \ldots, l_k \} \subset \{1, \ldots, n\}$ the notation $\MathVector{L}_{\mathcal{A} }$ means the vector $(L_{i_1}, \ldots, L_{i_k})$ while $\MathVector{L}_{\hat{\mathcal{A} } } = \MathVector{L}_{ \{1, \ldots, n\} \setminus \mathcal{A} }$.

There is one additional subtlety that crops up at this point. Note that for the case of $\Moduli_{1,1}(L)$, there is an order two automorphism obtained by rotating around the boundary by half a turn. There are two ways to deal with this issue. The approach taken in \cite{art:MirzakhaniWP} is to divide the appropriate integrals by 2 every time such a term appears in the above integral. Our approach, which is computationally equivalent, is to \emph{define} the volume of $\Moduli_{1,1}(L)$ to be half the value obtained by calculations using the above techniques. In other words, we have initial conditions
\begin{align*}
	\Vol_{0,3}(\MathVector{L}) &= 1 \\
	\Vol_{1,1}(L) &= \frac{1}{48}(L^2 + 4\pi^2).
\end{align*}
We will see that this viewpoint simplifies further calculations; as well, it agrees with known results from algebraic geometry.

The final step is to differentiate both sides with respect to $L_1$ and then integrate, which has the effect of simplifying the integrands on the right side of the equation. The result is
\begin{align*}
	\Vol_{g,n}(\MathVector{L}) 
	&= \frac{1}{2 L_1} \sum_{\substack{g_1+g_2 = g \\ \mathcal{A} \coprod \mathcal{B} } }
		\int_{0}^{L_1} \int_{0}^{\infty} \int_{0}^{\infty} xy H(t, x+y) \\[-2mm]
		& \hspace{45mm}\times \Vol_{g_1,n_1 }(x, \MathVector{L}_{\mathcal{A} } ) 
		\Vol_{g_2,n_2 }(y, \MathVector{L}_{\mathcal{B} }) dxdydt
	\\[2mm]
	&+\frac{1}{2L_1} \int_{0}^{L_1} \int_{0}^{\infty} \int_{0}^{\infty} xy H(t, x+y) \\
		&\hspace{35mm}\times \Vol_{g-1,n+1}(x,y,\MathVector{L}_{\hat{1} } )dxdydt
	\\
	&+\frac{1}{2L_1} \sum_{j=2}^{n} \int_{0}^{L_1}\int_{0}^{\infty} 
		x \bigl( H(x, L_1+L_j) + H(x, L_1 - L_j) \bigr) \\
		&\hspace{35mm}\times \Vol_{g,n-1}(x, \MathVector{L}_{\widehat{1,j } }) dxdt,
\end{align*}
where
\begin{equation*}
	H(x,y) = \frac{1 }{1 + e^{(x+y)/2 } }  +  \frac{1 }{1 + e^{(x-y)/2 } } .
\end{equation*}

\subsection{Relation to intersection numbers}

 In this subsection we review the idea of 
 Mirzakhani to write the integral over $\Moduli_{g,n}(\MathVector{L})$ as an integral of an appropriately modified volume form over $\Moduli_{g,n}$. This will relate the Weil-Petersson volumes to intersection numbers of tautological classes.
Following \cite{art:MirzakhaniIT}, let
\begin{equation*}
	\widehat{\Moduli}_{g,n} = \{(X, p_1, \ldots , p_n)\, |\, 
		X\in\CompactModuli_{g,n}(\MathVector{L}), \MathVector{L}\in\Reals_{\geq 0}^{n},
		p_i\in \beta_i \}
\end{equation*}
be the moduli space of bordered hyperbolic surfaces of arbitrary boundary length, with the additional information of a marked point on each boundary component. If $L_i = 0$, then we can think of $p_i$ as a point on a horocycle about the cusp.

\begin{figure}
	\includegraphics{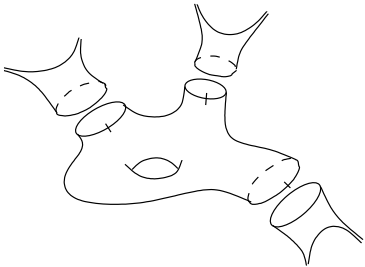}
	\caption{Capping off a bordered surface}
	\label{fig:SurfaceCap}
\end{figure}

The marked point can be used as a twist parameter, so by gluing on pairs of pants with two cusps and the third boundary having length matching the surface's boundary, we obtain a map $\widehat{\Moduli}_{g,n}\rightarrow \CompactModuli_{g,2n}$. In fact, we have $\widehat{\Moduli}_{g,n} = \CompactModuli_{g,2n}^{\Gamma}$ where $\Gamma = \{\gamma_1, \ldots, \gamma_n\}$ is a collection of curves which group the cusps into pairs. We refer to Figure~\ref{fig:SurfaceCap} for a descriptive picture of this construction.

This tells us that $\widehat{\Moduli}_{g,n}$ has a symplectic structure from the Weil-Petersson form on $\Moduli_{g,2n}^{\Gamma}$. Moreover it has a hamiltonian torus action given by rotating the marked points on the boundary. However, we need to take some care here. We are interested in studying the symplectic action in a neighborhood of  surfaces $X\in\widehat{\Moduli}_{g,n}$ with $l(\beta_i) = 0$. Moreover, we want to construct the action in such a way that these points are not fixed by the torus action. In other words, we need to non-trivially extend the action to the cusped surfaces. 

It is a simple matter of defining the twists to be proportional to the lengths of the boundaries. In other words, we scale the action so that a twist parameter of 1 is always the identity. The model is the change from the cartesian $(x,y)$ coordinates in the plane to the polar coordinate $(r,\theta)$. In the first case rotation around the origin leaves it fixed, but $(0,\theta)$ is not fixed by $\theta \mapsto\theta + \epsilon$.
From the point of view of $\widehat{\Moduli}_{g,n}$, the marking on the boundary degenerates to a marking on a horocycle of the cusp. The result, after a change of coordinates to the reparametrized twist coordinate $\theta_i = \tau_i/l_i$, is
\begin{equation*}
	\omega_{WP} = \sum l_i dl_i \wedge d\theta_i,
\end{equation*}
and the moment map corresponding to the twist vector field $\frac{\partial}{\partial\theta_i}$ is $\frac{1}{2}l_{i}^{2}$.

Given the map $\MathVector{L} : \widehat{\Moduli}_{g,n}  \rightarrow \Reals_{\geq 0}^{n}$ determined by mapping the marked surface $X\in\widehat{\Moduli}_{g,n}$ to the lengths of its boundary components, we see that $\MathVector{L}^{-1}(\MathVector{x})$ is a principal torus bundle over $\CompactModuli_{g,n}(\MathVector{x})$. In fact, over $\CompactModuli_{g,n}(0)$ it is the principal bundle associated to the vector bundle $\mathcal{L}_1 \oplus \cdots \oplus \mathcal{L}_n$.  At first glance this is a rather counter-intuitive statement, as marked points on the boundary map naturally to the tangent bundle, rather than the cotangent bundle. However, the principal torus bundle in question is naturally oriented from the induced orientations on the boundaries coming from the orientation of the surface. This orientation is opposite to the natural complex orientation on the tangent bundle. This is most easily seen by studying the clockwise orientation induced on the unit circle from the standard orientation of the plane.

Using symplectic reduction, we see that the reduced space $\MathVector{L}^{-1}(\MathVector{x} ) / T^n$ is symplectomorphic to $\CompactModuli_{g,n}(\MathVector{x})$ with the Weil-Petersson form. We may use the techniques of the Duistermaat--Heckman theorem to compare $\omega_{WP}(\MathVector{L})$ to $\omega_{WP}(0)$. The result is
\begin{equation*}
	\omega_{WP}(\MathVector{L}) = \omega_{WP}(0) - \frac{1}{2} \sum L^2_i\Curvature(\mathcal{L}_{i} ),
\end{equation*}
where $\Curvature(\mathcal{L}_{i})$ is the curvature of the bundle. Since $c_1(\mathcal{L}_i) = -\Curvature(\mathcal{L}_{i})$ we get
\begin{equation*}
	\omega_{WP}(\MathVector{L}) = \omega_{WP}(0) + \frac{1}{2} \sum L_{i}^{2} \psi_i.
\end{equation*}

\subsection{A rational recursion relation}

Using Wolpert's equivalence $\kappa_1 = \frac{\omega_{WP}}{2\pi^2}$ we define the rational volume of $\Moduli_{g,n}(\MathVector{L})$ to be
\begin{align*}
	v_{g,n}(\MathVector{L}) &\stackrel{\text{def}}{=} 
		\frac{\Vol_{g,n}(2\pi\MathVector{L})}{2^d\pi^{2d} } \quad\quad d = 3g-3+n \\
	&= \frac{1}{d!} \int_{\Moduli_{g,n}} ( \kappa_1 + \sum L_i^2 \psi_i )^d \\
	&= \sum_{\substack{d_0 + \cdots + d_n\\ = d}} \prod_{i=0}^{n}\frac{1}{d_i! } 
		\intersect{\kappa_1^{d_0} \prod\tau_{d_i} }_{g,n} \prod_{i=1}^{\infty}
		L_{i}^{2d_i}.
\end{align*}
We reformulate Mirzakhani's recursion relation for $\Vol_{g,n}$ into a recursion relation for $v_{g,n}$. Making the above change of variables to the recursion relation gives
\begin{equation*}
  \begin{split}
	v_{g,n}(\MathVector{L}) &=  \frac{2}{L_1} \int_{0}^{L_1} \int_{0}^{\infty} \int_{0}^{\infty}
		xy K(x+y, t) v_{g-1,n+1}(x,y,\MathVector{L}_{\hat{1}}) dxdydt \\
	&+ \frac{2}{L_1} \sum_{\substack{g_1+g_2=g\\ \mathcal{I} \coprod \mathcal{J} = \{2,\ldots,n\}} }
		\int_{0}^{L_1}\int_{0}^{\infty}\int_{0}^{\infty} xy K(x+y, t) \\[-3mm]
		&\hspace{51mm}\times v_{g_1,n_1}(x, \MathVector{L}_{\mathcal{I}}) 
			v_{g_2,n_2}(y, \MathVector{L}_{\mathcal{J}})dxdydt \\[2mm]
	&+ \frac{1}{L_1} \sum_{j=2}^{n} \int_{0}^{L_1} \int_{0}^{\infty} x 
		\left( K(x, t+L_j) + K(x, t-L_j)\right) \\[-2mm]
		&\hspace{30mm}\times v_{g, n-1}(x, \MathVector{L}_{\{\widehat{1,j}\} }) dxdt,
  \end{split}
\end{equation*}
with normalizations $v_{0,3}(L) = 1$ and $v_{1,1}(L) = \frac{1}{24}(1+L^2)$.

The integral kernel $K(x,t)$ is defined as
\begin{align}
	K(x,t) &= \frac{1}{1 + e^{\pi(x+t)}} + \frac{1}{1 + e^{\pi(x-t)}}, \notag \\
\intertext{which gives the following integral identities:}
	\hFunction_{2k+1}(t) &\stackrel{\text{def}}{=} \int_{0}^{\infty} \frac{x^{2k+1}}{(2k+1)!}
		K(x,t)dx \label{eqn:KernelSingleIntegral} \\
	&= \sum_{m=0}^{k+1} (-1)^{m-1} (2^{2m} - 2) \frac{B_{2m}}{(2m)!} \frac{t^{2k+2 - 2m}}{(2k+2-2m)!}, \notag \\
	\hFunction_{2i+2j+3}(t) &= \int_{0}^{\infty}\int_{0}^{\infty} \frac{x^{2i+1}y^{2j+1}}{(2i+1)! (2j+1)!}
		K(x+y,t)dxdy. \label{eqn:KernelDoubleIntegral}
\end{align}
Note that $B_{2m}$ is the $2m$-th Bernoulli number.

\section{From Mirzakhani's Recursion Relation to the Virasoro Algebra}
\label{sect:Virasoro}

Our aim is to show that the Mirzakhani recursion relations are equivalent to an algebraic constraint on the generating function for $\kappa_1$ and $\psi$ class intersections. We introduce the formal generating function for all $\kappa_1$ and $\psi$ class intersections
\begin{align*}
	G(s, t_0, t_1, t_2, \ldots) &\stackrel{\text{def}}{=}
		\sum_{g} \intersect{e^{s\kappa_1 + \sum t_i \tau_i} }_{g} \\
		&= \sum_{g} \sum_{m, \{n_i\}} \intersect{\kappa_1^m \tau_0^{n_0}\tau_1^{n_1}
			\cdots }_{g} 
		\frac{s^m}{m!} \prod_{i=0}^{\infty} \frac{t_{i}^{n_i}}{n_i! }.
\end{align*}
The main result of the paper is the following.
\begin{theorem}
	There exist a sequence of differential operators $V_{-1}, V_{0}, V_{1}, \ldots$ satisfying Virasoro relations
\begin{equation*}
	[V_{n}, V_{m}] = (n-m)V_{n+m}
\end{equation*}
and annihilating $\exp(G)$:
\begin{equation*}
	V_k \exp(G) = 0 \quad \text{for $k=-1, 0, 1, \ldots$}
\end{equation*}
This property uniquely fixes $G$ and enables one to calculate all coefficients of the expansion.
\end{theorem}

The proof is obtained by differentiating
 Mirzakhani's recursion relation. For reference, we note that
\begin{multline}
	\label{eqn:VolDerivatives}
	\frac{\partial^{2k_1}}{\partial L_1^{2k_1}} \cdots \frac{\partial^{2k_n}}{\partial L_n^{2k_n}}
		v_{g,n}(\MathVector{L}) \\
	= \sum_{\substack{d_0+\cdots + d_n = d \\ d_i \geq k_i}}  
		\frac{1}{d_0!} \prod_{i=1}^{n} 
		\left( \frac{(2d_i)!}{d_i! (2(d_i-k_i))!} L_{i}^{2(d_i - k_i)} \right) 
		\intersect{\kappa_1^{d_0} \tau_{d_1} \cdots \tau_{d_n}}_{g}, 
\end{multline}
and
\begin{equation}
	\label{eqn:VolDerivatives0}
	\frac{\partial^{2k_1}}{\partial L_1^{2k_1}} \cdots \frac{\partial^{2k_n}}{\partial L_n^{2k_n}}
		v_{g,n}(0) 
	= \frac{1}{k_0!} \prod_{i=1}^{n} \left( \frac{(2k_i)!}{k_i !} \right) \intersect{\kappa_1^{k_0} \tau_{k_1} \cdots \tau_{k_n}},
\end{equation}
where $k_0 = 3g - 3 + n - \sum_{i=1}^{n} k_i$.

The recursion relation gives the following identity for $(g,n) \neq (0,3), (1,1)$.

\begin{multline*}
	\frac{\partial^{2k_1}}{\partial L_1^{2k_1}} 
		\cdots \frac{\partial^{2k_n}}{\partial L_n^{2k_n}}
		v_{g,n}(0) \\
	\shoveleft{= \frac{\partial^{2k_1}}{\partial L_1^{2k_1}} \frac{2}{L_1} 
		\int_{0}^{L_1}\int_{0}^{\infty}\int_{0}^{\infty}
		xy K(x+y, t)} \\
			\shoveright{  \times 
				\frac{\partial^{2k_2}}{\partial L_2^{2k_2}} 
				\cdots \frac{\partial^{2k_n}}{\partial L_n^{2k_n}}
				v_{g-1,n+1}(x,y,\MathVector{L}_{\hat{1}}) dxdydt 
				\Bigr|_{\MathVector{L}=0}  } \\
	\shoveleft{  + \frac{\partial^{2k_1}}{\partial L_1^{2k_1}} \frac{2}{L_1} 
		\sum_{\substack{g_1+g_2=g \\ \mathcal{I}\coprod\mathcal{J}}}
		\int_{0}^{L_1}\int_{0}^{\infty}\int_{0}^{\infty} xy K(x+y, t) }\\
			\shoveright{   \times
				\frac{\partial^{2k(\mathcal{I})}}{\partial
				\MathVector{L}_{\mathcal{I}}^{2k(\mathcal{I})}}
				v_{g_1,n_1}(x, \MathVector{L}_{\mathcal{I}})
				\frac{\partial^{2k(\mathcal{J})}}{\partial 
				\MathVector{L}_{\mathcal{J}}^{2k(\mathcal{J})}}
				v_{g_2,n_2}(x, \MathVector{L}_{\mathcal{J}})dxdydt
				\Bigr|_{\MathVector{L}=0}  }\\
	\shoveleft{+ \sum_{j=2}^{n}
		\frac{\partial^{2(k_1+k_j)}}{\partial L_1^{2k_1}\partial L_j^{2k_j}}
		\frac{1}{L_1}  \int_{0}^{L_1} \int_{0}^{\infty} x( K(x, t+L_j) + K(x, t-L_j))  } \\
			\shoveright{  \times\frac{\partial^{2k(\widehat{1,j})}}{\partial 
				\MathVector{L}_{\widehat{1,j}}^{2k(\widehat{1,j})}}
				v_{g,n-1}(x, \MathVector{L}_{\widehat{1,j}}) dxdt
				\Bigr|_{\MathVector{L}=0}. } \\
\end{multline*}
Plugging in the expressions for derivatives of volume functions (\ref{eqn:VolDerivatives}), (\ref{eqn:VolDerivatives0}) and integrating against the kernel function using (\ref{eqn:KernelSingleIntegral}) and  (\ref{eqn:KernelDoubleIntegral}) gives
\begin{multline*}
	\frac{1}{k_0!}\prod_{i=1}^{n} \frac{(2k_i)!}{k_i!} \intersect{\kappa_1^{k_0}\tau_{k_1}\cdots \tau_{k_n}} \\
	\shoveleft{= \sum_{\stackrel{d_0 + d_1 + d_2}{= k_0 + k_1 - 2}}
		\frac{(2d_1 + 1)!(2d_2 +1)!}{d_0!d_1!d_2!} \prod_{i=2}^{n} \frac{(2k_i)!}{k_i!}
		\intersect{\kappa_1^{d_0}\tau_{d_1}\tau_{d_2}\MathVector{\tau}_{k(\hat{1})}}_{g-1,n+1} }\\
	\shoveright{   \times 
		\frac{\partial^{2k_1}}{\partial L_1^{2k_1}} \frac{2}{L_1}\int_{0}^{L_1} 
		\hFunction_{2(d_1 + d_2)+3}(t) dt \Bigr|_{L_1=0} }\\
	\shoveleft{   +\sum_{ \substack{ g_1+g_2=g \\ \mathcal{I}\coprod\mathcal{J} } }
		\sum_{\substack{ d_0+d_1=3g_1-3\\ +n_1 - k(\mathcal{I})\\
			d_{0}^{'} + d_{1}^{'}=3g_2-3 \\ +n_2 - k(\mathcal{J})  } }
		\frac{(2d_1+1)!(2d_{1}^{'} + 1)! }{d_0! d_1! d_{0}^{'}! d_{1}^{'}! }
		\prod_{i=2}^{n} \frac{(2k_i)! }{k_i!} 
		\intersect{\kappa_1^{d_0}\tau_{d_1}
			\MathVector{\tau}_{k(\mathcal{I})} }_{g_1,n_1} } \\
		\shoveright{  \times 
			\intersect{\kappa_1^{d^{'}_{0}}\tau_{d^{'}_{1}}
		\MathVector{\tau}_{k(\mathcal{J})} }_{g_2,n_2} 
		  \frac{\partial^{2k_1}}{\partial L_1^{2k_1}} \frac{2}{L_1}\int_{0}^{L_1} 
		\hFunction_{2(d_1 + d^{'}_{1})+3}(t) dt \Bigr|_{L_1=0} }\\
	\shoveleft{  + \sum_{j=2}^{n} \sum_{\substack{d_0 + d_1 = \\ k_0+k_1+k_j-1}}
		\frac{(2d_1+1)! }{d_0! d_1! } \prod_{i\neq 1,j} \frac{(2k_i)!}{k_i!}
		\intersect{\kappa_1^{d_0}\tau_{d_1}
		\MathVector{\tau}_{k(\widehat{1,j})} }_{g,n-1} }\\
		\shoveright{   \times \frac{\partial^{2k_1}}{\partial L_1^{2k_1} } \frac{2}{L_1}
		\int_{0}^{L_1} \hFunction^{(2k_j)}_{2d_1+1}(t) dt  \Bigr|_{L_1=0}.  }\\
\end{multline*}

We rewrite this sum by introducing the sequence of nonnegative integers $\{ n_0$, $n_1$, $n_2$, $\ldots\}$ such that $n_j = \FiniteCount{ \{k_i\, | \, i \neq 1, k_i=j\} }$, and relabel $k_1 = k$. In other words, we have
\begin{equation*}
	\intersect{\kappa_1^{k_0} \tau_{k_1} \cdots \tau_{k_n}}_{g} = 
		\intersect{\kappa_1^{k_0}\tau_k \tau_{0}^{n_0} \tau_{1}^{n_1} \cdots}_{g}.
\end{equation*}
We further define
\begin{equation}
	\beta_{i} = (-1)^{i-1} 2^i (2^{2i} - 2) \frac{B_{2i}}{(2i)!},
	\label{eqn:Beta}
\end{equation}
which results in the equation
\begin{multline}
	(2k+1)!! \intersect{\kappa_1^{k_0}\tau_{k} {\prod_{i=0}^{\infty}}\tau_{i}^{n_i} }_{g} \\
	= \frac{1}{2} \sum_{\substack{d_0+d_1+d_2 =\\ k_0+k-2} }\frac{k_0!}{d_0!}
		\beta_{(k_0-d_0)} (2d_1 +1)! (2d_2+1)! 
		\intersect{\kappa_1^{d_0}\tau_{d_1}\tau_{d_2} 
			{\prod_{i=0}^{\infty}}\tau_{i}^{n_i} }_{g-1,n+1} \\
	+ \frac{1}{2} \sum_{ \substack{g_1+g_2=g \\ \{l_i\} + \{m_i\} = \{n_i\} } }
		\sum_{\substack{ d_0+d_1=3g_1-3\\ +n_1 - k(\mathcal{I})\\
			d_{0}^{'} + d_{1}^{'}=3g_2-3 \\ +n_2 - k(\mathcal{J})  } }
		\frac{k_0! }{d_0!d^{'}_{0}! } \beta_{(k_0-d_0-d^{'}_{0})}
		(2d_1 + 1)!(2d^{'}_{1} + 1)! \\
	\shoveright{ \times
		\prod_{i=0}^{\infty} \frac{n_i! }{l_i! m_i! }
		\intersect{\kappa_1^{d_0}\tau_{d_1} \prod \tau_{i}^{l_i} }_{g_1}
		\intersect{\kappa_1^{d^{'}_{0}}\tau_{d^{'}_{1}} \prod \tau_{i}^{m_i} }_{g_2} }\\
	\shoveleft{ + \sum_{j=0}^{\infty} \sum_{\substack{ d_0+d_1 = \\ k_0+k+j-1 } }
		\frac{k_0!}{d_0!} \beta_{(k_0-d_0) } \frac{(2d_1+1)! }{(2j-1)! }n_j
		\intersect{\kappa_1^{d_0}\tau_{d_1}\tau^{-1}_{j}\prod \tau_{i}^{n_i} }_{g}
		\label{eqn:RecursionExpression}. } \\
\end{multline}

By looking at expressions of the form
\begin{align*}
	\frac{\partial G}{\partial t_i} &= \sum_{g} \sum_{m, \{n_i\}} \intersect{\kappa_1^m \tau_i \tau_0^{n_0}\tau_1^{n_1}
			\cdots }_{g} 
		\frac{s^m}{m!} \prod_{i=0}^{\infty} \frac{t_{i}^{n_i}}{n_i! }, \\
	s^i t_j \frac{\partial G}{\partial t_k} &= \sum_{g} \sum_{m, \{n_i\}} \frac{m!}{(m-i)!}n_j \intersect{\kappa_1^{m-i} \tau_j^{-1} \tau_k \tau_0^{n_0}\tau_1^{n_1}
			\cdots }_{g} 
		\frac{s^m}{m!} \prod_{i=0}^{\infty} \frac{t_{i}^{n_i}}{n_i! },  \\
	s^i \frac{\partial^2 G }{\partial t_j \partial t_k } &= \sum_{g} \sum_{m, \{n_i\}}
	\frac{m! }{(m-i)! }
	\intersect{\kappa_1^{m-i} \tau_j \tau_k \tau_0^{n_0}\tau_1^{n_1} \cdots }_{g} 
		\frac{s^m}{m!} \prod_{i=0}^{\infty} \frac{t_{i}^{n_i}}{n_i! }, \\
	s^i \frac{\partial G }{\partial t_j} \frac{\partial G }{\partial t_k } &=
		\sum_{\substack{g \\ m, \{n_i\} } } 
		\sum_{\substack{g_1 + g_2 = g \\ d_1 + d_2 = m-i \\ \{k_i\} + \{l_i\} = \{n_i\} } }
		\frac{m! }{d_1! d_2! }
		 \intersect{\kappa_1^{d_1} \tau_j \tau_0^{k_0} \cdots }_{g_1} 
		 \intersect{\kappa_1^{d_2} \tau_k \tau_0^{l_0} \cdots }_{g_2}
		\frac{s^m}{m!} \prod_{i=0}^{\infty} \frac{t_{i}^{n_i}}{n_i! },
\end{align*}
we see that (\ref{eqn:RecursionExpression}) leads to the following expression for all $k>0$: 
\begin{multline*}
	(2k+3)!! \frac{\partial}{\partial t_{k+1}}
		G = \sum_{i,j=0}^{\infty} 
		\frac{(2(i+j+k) +1)!! }{(2j-1)!! } \beta_{i} s^i t_j \frac{\partial}{\partial t_{i+j+k}}G \\
	 + \frac{1}{2} \sum_{i=0}^{\infty} \sum_{\substack{d_1+d_2= \\i+k-1}} (2d_1+1)!!
	 	(2d_2+1)!! \beta_{i}s^i \left( \frac{\partial^2 G}{\partial t_{d_1} \partial t_{d_2} }
		+ \frac{\partial G }{\partial t_{d_1} }\frac{\partial G }{\partial t_{d_2} }  \right).
\end{multline*}
Note that  similar expressions are possible for $k=-1,0$ by taking special care of the base cases $(g,n) = (0,3), (1,1)$. We introduce the family of differential operators for $k \geq -1$
\begin{multline*}
 \DeformedVirasoro_{k} = -\frac{(2k+3)!! }{2}\frac{\partial}{\partial t_{k+1}}
 	+ \delta_{k,-1}(\frac{t_0^2}{4} + \frac{s}{48}) + \frac{\delta_{k,0}}{48} \\
 	+ \frac{1}{2}\sum_{i,j=0}^{\infty} \frac{(2(i+j+k)+1)!! }{(2j-1)!! }\beta_{i}s^i t_j 
	\frac{\partial}{\partial t_{i+j+k}} \\
	+ \frac{1}{4}\sum_{i=0}^{\infty} \sum_{\substack{d_1+d_2=\\i+k-1 }}(2d_1 +1)!!
		(2d_2+1)!! \beta_{i} s^i \frac{\partial^2}{\partial t_{d_1} \partial t_{d_2}  }.
\end{multline*}
We have proven the following statement.

\begin{theorem} For $k \geq -1$
  \begin{equation*}
	\DeformedVirasoro_k \exp(G) = 0.
  \end{equation*}
\end{theorem}

A reasonable question is: what is the algebra spanned by the operators $\DeformedVirasoro_k$? The answer is that they span a subalgebra of the Virasoro algebra. One can check directly that the operators satisfy the relations
\begin{equation*}
	[\DeformedVirasoro_{n}, \DeformedVirasoro_{m}] = 
		(n-m)\sum_{i=0}^{\infty} \beta_{i} s^{i} \DeformedVirasoro_{n+m+i}.
\end{equation*}
On the surface, this looks to be a deformation of the Virasoro relations (setting $s=0$ recovers Virasoro). However, this is, in fact, a simple reparametrization of the representation. These statements can all be proved by direct calculations. Here we make some simplifications.
Let us introduce new variables $\{T_{2j+1}\}_{j=0,1,\ldots}$ defined by
\begin{equation*}
	T_{2i+1} = \frac{t_i}{(2i+1)!!},
\end{equation*}
which transform the operators $\DeformedVirasoro_k$ into
\begin{align*}
	\DeformedVirasoro_{k} &= -\frac{1}{2}\frac{\partial }{\partial T_{2k+3} }
	+ \delta_{k,-1} (\frac{t_0^2}{4} + \frac{s}{48}) + \frac{\delta_{k,0} }{16 } \\
	&\quad + \frac{1}{2} \sum_{i,j=0}^{\infty} (2j+1) \beta_{i} s^i T_{2j+1} \frac{\partial}{\partial 	T_{2(i+j+k) + 1}} \\
	&\quad + \frac{1}{4} \sum_{i=0}^{\infty} \sum_{\substack{d_1+d_2 = \\ i+k-1 }} \beta_{i} s^i \frac{\partial^2}{\partial T_{2d_1 +1} \partial T_{2d_2 +1} }.
\end{align*}
This admits the following `boson' representation, similar to that used by Kac and Schwarz \cite{art:KacSchwarz}. Define operators $J_p$ for $p\in\Integers$ by
\begin{equation*}
	J_p = \begin{cases}
		(-p) T_{-p} & \text{if $p<0$}, \\
		\frac{\partial}{\partial T_p} & \text{if $p>0$}.
	\end{cases}
\end{equation*}
Then
\begin{align*}
	\DeformedVirasoro_k &= -\frac{1}{2} J_{2k+3} + \sum_{i=0}^{\infty}
		\beta_{i}s^i E_{k+i}, \\
	\intertext{where}
	E_k &= \frac{1}{4}\sum_{p\in\Integers} J_{2p+1} J_{2(k-p) - 1} 
		+ \frac{\delta_{k,0}}{16}.
\end{align*}
To recover operators satisfying the Virasoro constraint we need a better handle on the constants $\beta_{i}$, as defined by (\ref{eqn:Beta}).  Starting from the defining formula for the Bernoulli numbers
\begin{equation*}
	\sum_{n=0}^{\infty} \frac{B_{2n}}{(2n)! } z^{2n}
	= \frac{z}{2}\frac{e^{z/2} + e^{-z/2} }{e^{z/2} - e^{-z/2} },
\end{equation*}	
we see that
\begin{align*}
	\sum_{i=0}^{\infty} \beta_{i}s^i &= \sqrt{2s} (\cot \sqrt{s/2} - \cot \sqrt{2s}) \\
	&= \frac{\sqrt{2s} }{\sin\sqrt{2s} }.
\end{align*}
This motivates the definition of the constants $\alpha_{i}$ by the series
\begin{equation*}
	\sum_{i=0}^{\infty} \alpha_{i} s^i = \frac{\sin\sqrt{2s} }{\sqrt{2s} },
\end{equation*}
from which we obtain the operators
\begin{align}
	\Virasoro_k &\stackrel{\text{def}}{=} \sum_{i=0}^{\infty} 
		\alpha_{i}s^i \DeformedVirasoro_{k+i} \notag \\
	&= -\frac{1}{2}\sum_{i=0}^{\infty} \alpha_{i}s^i J_{2k+3} + E_k. \label{eqn:Virasoro}
\end{align}
We are now ready to prove the following.
\begin{proposition}
	The operators $\Virasoro_k$, $k\geq -1$ satisfy the Virasoro relations
\begin{equation*}
	[\Virasoro_n, \Virasoro_m] = (n-m)\Virasoro_{n+m}.
\end{equation*}
\begin{proof}
	The first step is to verify that operators $E_k$ satisfy the Virasoro relations, which is a straightforward calculation. Since	$[J_{2k+3}, E_m] =  (2k+3)J_{2(k+m)+3}$ we see that
\begin{align*}
	[\Virasoro_n, \Virasoro_m] &= \biggl[-\frac{1}{2} \sum_{i=0}^{\infty}\alpha_{i} s^i
		J_{2(n+i)+3 } + E_n, -\frac{1}{2}\sum_{j=0 }^{\infty} \alpha_{i} s^i
		J_{2(m+j)+3 } + E_m\biggr] \\
		&= -\frac{1}{2} \sum_{i=0}^{\infty} \alpha_{i} s^i \biggl( \bigl[
			J_{2(n+i)+3 }, E_m
		\bigr] + \bigl[ 
			E_n, J_{2(m+i)+3 }
		\bigr] \biggr) \\
		&= (n-m) \Virasoro_{n+m }.
\end{align*}
\end{proof}
\end{proposition}

 \section{Relationship to KdV hierarchy}
 \label{sect:IntegrableSystem}
 
The Witten-Kontsevich theorem  \cite{art:Witten2dGravity, art:Kontsevich}
states that the generating function for $\psi$ class intersections
\begin{align*}
	F(t_0, t_1, \ldots) &= \sum_{g} \intersect{e^{\sum \tau_i t_i} }_{g} \\
		&= \sum_{g} \sum_{\{n_{*} \} } \intersect{\prod \tau_{i}^{n_i} }_{g} \prod \frac{t_i^{n_i} }{n_i! }
\end{align*}
is a $\tau$-function for the KdV hierarchy.
The property of being a $\tau$ function, combined with the string equation
 \begin{equation*}
 	\intersect{\tau_0 \prod_{i=1}^{n} \tau_{d_i} }_{g} = \sum_{j=1}^{n} 
		\intersect{\prod_{i=1}^{n} \tau_{d_i - \delta_{ij}} }_{g},
 \end{equation*}
 completely determines the function $F$. Another way of 
 determining $F$ is the Virasoro
 constraint condition. Let us
  define the sequence of operators $L_k$ for $k \geq -1$:
 \begin{multline}
 	\label{eqn:WK-Virasoro}
 	L_{k} = -\frac{(2k+3)!! }{2}\frac{\partial }{\partial t_{k+1} } 
				+ \frac{1}{2} \sum_{j=0}^{\infty} \frac{(2(j+k) + 1)!! }{(2j-1)!! } t_j 
		\frac{\partial }{\partial t_{j+k}} \\
	+ \frac{1}{4} \sum_{\substack{d_1 + d_2 = k-1 \\ d_1,d_2 \geq 0}} 
		(2d_1 + 1)!!(2d_2 + 1)!! 
		\frac{\partial^2 }{\partial t_{d_1} \partial t_{d_2} }
		+ \frac{\delta_{k,-1} t_0^2}{4} + \frac{\delta_{k,0} }{48}. 
 \end{multline}
 The Witten-Kontesevich theorem, together with the string equation, implies 
 	$$
	L_{k} (\exp F) = 0
	$$ for $k \geq -1$.
This property is also sufficient to uniquely fix $F$. Note that $L_{-1} e^F = 0 $ is equivalent to the string equation. The consistency of the infinite set of differential equations follows from the fact that operators $L_n$ satisfy the Virasoro relations:
 \begin{equation*}
 	[L_n, L_m] = (n-m)L_{n+m}.
 \end{equation*}
 Recall the operators $V_k$ defined in equation~\ref{eqn:Virasoro} (rewritten in terms of the variables $t_i$)
 \begin{multline*}
 	\Virasoro_k = -\frac{1}{2} \sum_{i=0}^{\infty} (2(i+k)+3)!! \alpha_{i} s^i 
		\frac{\partial }{\partial t_{i+k+1} } 
		+ \frac{1}{2} \sum_{j=0}^{\infty} \frac{(2(j+k)+1)!! }{(2j-1)!! } t_j
		\frac{\partial }{\partial t_{j+k} } \\
	+ \frac{1}{4} \sum_{\substack{d_1 + d_2 = k-1 \\ d_1, d_2 \geq 0 } }
		(2d_1 + 1)!! (2d_2 + 1)!! \frac{\partial^2 }{\partial t_{d_1} \partial t_{d_2}}
		+ \frac{\delta_{k,-1}t_0^2}{4} + \frac{\delta_{k,0} }{48}, 
 \end{multline*}
 where $\alpha_{i} = \frac{(-2)^i }{(2i+1)! }$.
The change of variables
 \begin{equation*}
 	\tilde{t}_i = \begin{cases}
		t_i & \text{for $i=0,1$ }, \\
		t_i - (2i-1)!! \alpha_{i-1} s^{i-1} & \text{otherwise,}
	\end{cases}
 \end{equation*}
transforms the operators $V_k$ into
\begin{multline*}
 	\Virasoro_k = -\frac{1}{2}  (2k+3)!! 
		\frac{\partial }{\partial \tilde{t}_{k+1} } 
		+ \frac{1}{2} \sum_{j=0}^{\infty} \frac{(2(j+k)+1)!! }{(2j-1)!! } \tilde{t}_j
		\frac{\partial }{\partial \tilde{t}_{j+k} } \\
	+ \frac{1}{4} \sum_{\substack{d_1 + d_2 = k-1 \\ d_1, d_2 \geq 0 } }
		(2d_1 + 1)!! (2d_2 + 1)!! \frac{\partial^2 }{\partial \tilde{t}_{d_1} 
		\partial \tilde{t}_{d_2}}
		+ \frac{\delta_{k,-1}\tilde{t}_0^2}{4} + \frac{\delta_{k,0} }{48}. 
 \end{multline*}
 But these are precisely the operators $L_k$ (\ref{eqn:WK-Virasoro}). We have thus proven the following.
 \begin{theorem}
 	\begin{equation*}
		G(s, t_0, t_1, \ldots) = F(t_0, t_1, t_2 + \gamma_2, t_3 + \gamma_3, \ldots),
	\end{equation*}
where $\gamma_i = \frac{(-1)^i}{(2i+1)i!}s^{i-1}$. In particular, for any fixed value of $s$, $G$ is a $\tau$ function for the KdV hierarchy.
\label{theorem:KdV}
 \end{theorem}

That the more general generating function $G$ is expressible in terms of $F$ is not a surprise. It has been known since at least the work of Witten \cite{art:Witten2dGravity} that intersections involving $\kappa$ classes are expressible in terms of $\psi$ classes. Moreover, Faber's formula \cite{art:Faber} for this correspondence gives an explicit proof of the above theorem. In fact, one has
\begin{equation*}
	\kappa_{1}^{n} = \sum_{\substack{\sigma \in S_n \\ (\sigma = \gamma_1\cdots 		\gamma_k \\ \text{is cycle decomp}) } }
	\frac{(-1)^{n-k} }{ {\textstyle \prod_{i=1}^{k}} (|\gamma_i | -1 )!}
	\pi_{\{q_1, \ldots q_k \}_* } (\psi_{q_1}^{|\gamma_1| + 1 } \cdots
	\psi_{q_k}^{|\gamma_k | + 1 } ),
\end{equation*}
which gives a short, direct proof of Theorem~\ref{theorem:KdV}. This is essentially the approach taken by Zograf \cite{art:Zograf} for his calculation of the Weil-Petersson volumes of $\Moduli_{g,n}$.

\bibliographystyle{hamsplain}
\bibliography{References}

\end{document}